\def\author{S L Kleiman}
\def\title{Equisingularity, Multiplicity, and Dependence}
\def\date{13 May 1998}
\def\TheMagstep{\magstep1}      
\def\PaperSize{letter}          
\def\abstract
 {This is a report on some recent work by Gaffney, Massey, and the
author, characterizing the conditions \Af\ and \Wf\ for a family of ICIS
germs equipped with a function.  First we introduce the work informally.
Then we review the formal definitions of \Af\ and \Wf, and state the
theorems that characterize them by the constancy of Milnor numbers.
Next we review the definition of the Buchsbaum--Rim multiplicity, and
reformulate the theorems by the constancy of certain Buchsbaum--Rim
multiplicities.  Finally, we review the theory of integral dependence of
elements on submodules of free modules, and apply it to prove the
reformulated theorems.}
\let\:=\colon  \let\x=\times 
\def\Wf{\hbox{\rm W$_f$}}	\def\Af{\hbox{\rm A$_f$}}
\def\Wfi{\hbox{\rm W$_{f^i}$}}	\def\Afi{\hbox{\rm A$_{f^i}$}}
 \let\?=\overline 
 
 \def\IP{{\bf P}} \def\IC{{\bf C}} \def\Im{{\bf m}}

\def\mcl#1{\expandafter\def\csname c#1\endcsname{{\cal#1}}}
 \mcl E \mcl F \mcl I \mcl J \mcl M \mcl N \mcl O
\def\mop#1 {\expandafter\def\csname#1\endcsname
  {\mathop{\rm#1}\nolimits}}
\mop Projan  \mop Supp  \mop dist

\def\and{\hbox{\rm\ and }} 
\def\pd#1/#2 {\partial#1/\partial#2}
\def\tc#1 #2 {{\textstyle{#1\choose#2}}}

 \let\@=@  

\parskip=0pt plus 1.75pt \parindent10pt
\hsize29pc 
\vsize43pc 
\abovedisplayskip 4pt plus3pt minus1pt
\belowdisplayskip=\abovedisplayskip
\abovedisplayshortskip 2.5pt plus2pt minus1pt
\belowdisplayshortskip=\abovedisplayskip

\def\TRUE{TRUE} 
\ifx\DoublepageOutput\TRUE \def\TheMagstep{\magstep0} \fi
\mag=\TheMagstep

\newskip\vadjustskip \vadjustskip=0.5\normalbaselineskip
\def\centertext
 {\hoffset=\pgwidth \advance\hoffset-\hsize
  \advance\hoffset-2truein \divide\hoffset by 2\relax
  \voffset=\pgheight \advance\voffset-\vsize
  \advance\voffset-2truein \divide\voffset by 2\relax
  \advance\voffset\vadjustskip
 }
\newdimen\pgwidth\newdimen\pgheight
\def\letter{letter}\def\AFour{AFour}
\ifx\PaperSize\letter
 \pgwidth=8.5truein \pgheight=11truein
 \message{- Got a paper size of letter.  }\centertext
\fi
\ifx\PaperSize\AFour
 \pgwidth=210truemm \pgheight=297truemm
 \message{- Got a paper size of AFour.  }\centertext
\fi

\def\today{\ifcase\month\or	
 January\or February\or March\or April\or May\or June\or
 July\or August\or September\or October\or November\or December\fi
 \space\number\day, \number\year}
\nopagenumbers
 \newcount\pagenumber \pagenumber=1
 \def\advancepagenumber{\global\advance\pagenumber by 1}
\def\folio{\number\pagenum} 
\headline={%
  \ifnum\pagenum=0\hfill
  \else
   \ifnum\pagenum=1\firstheadline
   \else
     \ifodd\pagenum\oddheadline
     \else\evenheadline\fi
   \fi
  \fi
}
\expandafter\ifx\csname date\endcsname\relax \let\dato=\today
	    \else\let\dato=\date\fi
\def\oddheadline{\eightpoint \rlap{\dato}
 \hfil\headtitle\hfil\llap{\folio}}
\def\evenheadline{\eightpoint\rlap{\folio}
 \hfil\author\hfil\llap{\dato}}
\def\headtitle{\title}

 \newdimen\fullhsize \newbox\leftcolumn
 \def\fulline{\hbox to \fullhsize}
\def\doublepageoutput
{\let\lr=L
 \output={\if L\lr
           \global\setbox\leftcolumn=\columnbox \global\let\lr=R%
          \else \doubleformat \global\let\lr=L
	  \fi
        \ifnum\outputpenalty>-20000 \else\dosupereject\fi
	}%
 \def\doubleformat{\shipout\vbox{%
     \ifx\PaperSize\AFour
	   \fulline{\hfil\box\leftcolumn\hfil\columnbox\hfil}%
     \else
	   \fulline{\hfil\hfil\box\leftcolumn\hfil\columnbox\hfil\hfil}%
     \fi             }%
     \advancepageno
}
 \def\columnbox{\vbox
   {\if E\topmark\headline={\hfil}\nopagenumbers\fi
    \makeheadline\pagebody\makefootline\advancepagenumber}%
   }%
\fullhsize=\pgheight \hoffset=-1truein
 \voffset=\pgwidth \advance\voffset-\vsize
  \advance\voffset-2truein \divide\voffset by 2
  \advance\voffset\vadjustskip
 
  \pagenum=0\null\vfill\nopagenumbers\eject\pagenum=1\relax
}
\ifx\DoublepageOutput\TRUE \let\pagenum=\pagenumber\doublepageoutput
 \else \let\pagenum=\pageno \fi

 \font\twelvebf=cmbx12          
 \font\smc=cmcsc10              
\catcode`\@=11          
\def\eightpoint{\eightpointfonts
 \setbox\strutbox\hbox{\vrule height7\p@ depth2\p@ width\z@}%
 \eightpointparameters\eightpointfamilies
 \normalbaselines\rm
 }
\def\eightpointparameters{%
 \normalbaselineskip9\p@
 \abovedisplayskip9\p@ plus2.4\p@ minus6.2\p@
 \belowdisplayskip9\p@ plus2.4\p@ minus6.2\p@
 \abovedisplayshortskip\z@ plus2.4\p@
 \belowdisplayshortskip5.6\p@ plus2.4\p@ minus3.2\p@
 }
\newfam\smcfam
\def\eightpointfonts{%
 \font\eightrm=cmr8 \font\sixrm=cmr6
 \font\eightbf=cmbx8 \font\sixbf=cmbx6
 \font\eightit=cmti8
 \font\eightsmc=cmcsc8
 \font\eighti=cmmi8 \font\sixi=cmmi6
 \font\eightsy=cmsy8 \font\sixsy=cmsy6
 \font\eightsl=cmsl8 \font\eighttt=cmtt8}
\def\eightpointfamilies{%
 \textfont\z@\eightrm \scriptfont\z@\sixrm  \scriptscriptfont\z@\fiverm
 \textfont\@ne\eighti \scriptfont\@ne\sixi  \scriptscriptfont\@ne\fivei
 \textfont\tw@\eightsy \scriptfont\tw@\sixsy \scriptscriptfont\tw@\fivesy
 \textfont\thr@@\tenex \scriptfont\thr@@\tenex\scriptscriptfont\thr@@\tenex
 \textfont\itfam\eightit        \def\it{\fam\itfam\eightit}%
 \textfont\slfam\eightsl        \def\sl{\fam\slfam\eightsl}%
 \textfont\ttfam\eighttt        \def\tt{\fam\ttfam\eighttt}%
 \textfont\smcfam\eightsmc      \def\smc{\fam\smcfam\eightsmc}%
 \textfont\bffam\eightbf \scriptfont\bffam\sixbf
   \scriptscriptfont\bffam\fivebf       \def\bf{\fam\bffam\eightbf}%
 \def\rm{\fam0\eightrm}%
 }
\def\vfootnote#1{\insert\footins\bgroup
 \eightpoint\catcode`\^^M=5\leftskip=0pt\rightskip=\leftskip
 \interlinepenalty\interfootnotelinepenalty
  \splittopskip\ht\strutbox 
  \splitmaxdepth\dp\strutbox \floatingpenalty\@MM
  \leftskip\z@skip \rightskip\z@skip \spaceskip\z@skip \xspaceskip\z@skip
  \textindent{#1}\footstrut\futurelet\next\fo@t}

\def\p.{p.\penalty\@M \thinspace}
\def\pp.{pp.\penalty\@M \thinspace}
\newcount\sctno
\def\sctn#1\par
  {\removelastskip\vskip0pt plus5\normalbaselineskip \penalty250
  \vskip0pt plus-5\normalbaselineskip \bigskip\bigskip
  \centerline{#1}\bigskip}
\def\sct#1.{\sctno=#1\relax\sctn#1.}

\def\item#1 {\par\indent\indent\indent
 \hangindent3\parindent
 \llap{\rm (#1)\enspace}\ignorespaces}
 \def\inpart#1 {{\rm (#1)\enspace}\ignorespaces}
 \def\part {\par\inpart}

\def\Cs#1){\(\number\sctno.#1)}
\def\proclaim#1 #2 {\medbreak
  {\bf#1 (\number\sctno.#2)}\enspace\bgroup
\it}
\def\endproclaim{\par\egroup\medskip}

\def\lem{\proclaim Lemma } \def\prp{\proclaim Proposition }
\def\cor{\proclaim Corollary }	\def\thm{\proclaim Theorem }

\def\dno#1${\eqno\hbox{\rm(\number\sctno.#1)}$}
\def\Cs#1){\unskip~{\rm(\number\sctno.#1)}}
\def\dleqno#1{\hfil\llap{\rm#1}\hfilneg}

 \newcount\refno \refno=0        \def\NoKey{*!*}
 \def\MakeKey{\advance\refno by 1 \expandafter\xdef
  \csname\TheKey\endcsname{{\number\refno}}\NextKey}
 \def\NextKey#1 {\def\TheKey{#1}\ifx\TheKey\NoKey\let\next\relax
  \else\let\next\MakeKey \fi \next}
 \def\RefKeys #1\endRefKeys{\expandafter\NextKey #1 *!* }
 \def\SetRef#1 #2,{\hang\llap
  {[\csname#1\endcsname]\enspace}{\smc #2},}
 \newbox\keybox \setbox\keybox=\hbox{[25]\enspace}
 \newdimen\keyindent \keyindent=\wd\keybox
\def\references{\kern-\medskipamount
  \sctn References\par
  \vskip-\medskipamount
  \bgroup   \frenchspacing   \eightpoint
   \parindent=\keyindent  \parskip=\smallskipamount
   \everypar={\SetRef}\par}
\def\endreferences{\egroup}

 \def\serial#1#2{\expandafter\def\csname#1\endcsname ##1 ##2 ##3
        {\unskip\ {\it #2\/} {\bf##1} (##2), ##3.}} 
 \serial{ajm}{Amer. J. Math.}
 \serial{asens}{Ann. Scient. \'Ec. Norm. Sup.}
 \serial{cmh}{Comment. Math. Helvetici}
 \serial{conm}{Contemp. Math.}
 \serial{crasp}{C. R. Acad. Sci. Paris}
 \serial{faa}{Funct. Anal. Appl.}
 \serial{invent}{Invent. Math.}
 \serial{ja}{J. Algebra}
 \serial{ma}{Math. Ann.}
 \serial{mpcps}{Math. Proc. Camb. Phil. Soc.}
 \serial{splm}{Springer Lecture Notes in Math.}
 \serial{tams}{Trans. Amer. Math. Soc.}

\def\UThin{\penalty\@M \thinspace\ignorespaces}
\def\(#1){{\let~=\UThin\rm(#1)}}
\def\relaxnext@{\let\next\relax}
\def\cite#1{\relaxnext@
 \def\nextiii@##1,##2\end@{\unskip\space{\rm[\SetKey{##1},\let~=\UThin##2]}}%
 \in@,{#1}\ifin@\def\next{\nextiii@#1\end@}\else
 \def\next{{\rm[\SetKey{#1}]}}\fi\next}
\newif\ifin@
\def\in@#1#2{\def\in@@##1#1##2##3\in@@
 {\ifx\in@##2\in@false\else\in@true\fi}%
 \in@@#2#1\in@\in@@}
\def\SetKey#1{{\bf\csname#1\endcsname}}

\catcode`\@=12  

 \RefKeys
 BMM BR63 Ga92 Ga96 GK98 GM98 GnM96 Gr75 HM82 HM94 HMS84 H76 KR94 KT94 KT96
KT97 Le74 LS73 LT88 Lo84 M98 N80 Re87 Te73 Te80 Te81 V76
 \endRefKeys
   
 \def\firstheadline{alg-geom/9805062\hfill
	{\sl In onore di Mario FIORENTINI}}
{\leftskip=0pt plus1fill \rightskip=\leftskip
 \obeylines
 \leavevmode \bigskip\bigskip\bigskip
 {\twelvebf \title
 } \bigskip
 \footnote{}{\noindent 
  Subj-class: 32S15 (Primary) 14B05, 13H15 (Secondary).\par
  Research supported in part by NSF grant 9400918-DMS.}
 Steven L Kleiman
   \footnote{}{%
 It is a great pleasure for the author to thank Terence Gaffney for
enthusiastically introducing him to the present theory, for patiently
answering his countless questions, and for generously encouraging their
fruitful collaboration.  It is also a pleasure to thank David Massey for
several helpful discussions.  In addition, the author is sincerely
grateful to Herwig Hauser and Bernard Teissier for reading an earlier
draft and making valuable comments and suggestions.}
  {\eightpoint\it\bigskip
 Department of Mathematics, Room {\sl 2-278} MIT,
 {\sl77} Mass Ave, Cambridge, MA {\sl02139-4307}, USA
 \rm E-mail: \tt Kleiman\@math.mit.edu \bigskip
 \rm \dato \bigskip
 }
}
{\parindent=1.5\parindent \narrower \eightpoint \noindent
 {\smc Abstract.}\enspace \abstract \par}


\sct1. Introduction

The conditions \Af\ and \Wf\ are ``relative'' forms of the Whitney
conditions A and B (or W).  The latter are important as they are
manageable algebraic-geometric conditions on the tangents to a variety
$X$, yet they imply the local topological triviality of $X$ along a
given smooth subvariety $Y$ such that $X-Y$ is smooth.  When $X$ carries
a map $f$ that has constant rank off $Y$, then \Af\ and \Wf\ are defined
as the corresponding conditions on the tangents to the level sets (or
fibers) of $f$.  Thom introduced \Af\ as the primary condition that
would ensure the topological triviality of the pair $X,f$ along $Y$.

By definition, \Af\ and \Wf\ signify this.  For simplicity, assume that
$f$ is a nonconstant function vanishing on $Y$, and embed $X$ in $\IC^n$
so that $Y$ is a linear subspace through the origin $0$.  Then \Af\
holds at $0$ if $Y$ lies in every hyperplane obtained as a limit of
hyperplanes $H$, each tangent to a level hypersurface of $f$ at a point
$x$ of $X-Y$, as $x$ aproaches $0$.  A more stringent condition, \Wf\
signifies that the angle between $H$ and $Y$ approaches $0$ as fast as
$x$ approaches $0$.

The Thom--Mather second isotopy lemma readily yields the following
isotopy theorem: if \Wf\ holds at $0$, then, after $X$ is replaced by a
neighborhood of $0$, the pair $X,f$ is topologically right trivial.
More precisely, let $X_y$ denote the fiber of a transverse projection to
$Y$, and set $f_y:=f|X_y$.  Then there is a homeomorphism $h$ from the
product $X_0\x Y$ onto $X$ such that $fh$ is equal to $f_0\x1_Y$.  Thus,
if \Wf\ holds at $0$, then, for $y\in Y$ near 0, the pairs $X_y,f_y$ are
topologically the same, or ``topologically equisingular.''

In the present article, we discuss and develop some recent work done by
Gaffney, Massey, and the author in \cite{GK98} and \cite{GM98} on the
algebraic-geometric significance of \Af\ and \Wf.  Here $X$ is the total
space of a complex analytic family of germs of isolated
complete-intersection singularities (ICIS germs) $X_y$, parameterized by
a smooth variety $Y$.

For convenience, we identify $Y$ with the subvariety of $X$ traced by
the central points $0\in X_y$ as $y$ varies.  We assume $f$ vanishes on
$Y$.  We choose an embedding of $X$ in $\IC^a\x\IC^b$ such that $Y$
represents the germ of $0\x\IC^b$.  We choose an extension of $f$ over a
neighborhood of $X$, and we denote the extension by $f$ too.  Needless
to say, our results are independent of these choices.

What does \Af\ mean by itself?  Let $Z_y$ be the level hypersurface
through 0 of $f_y$.  It turns out that \Af\ is closely related to the
vanishing cycles associated to the individual $X_y$ and $Z_y$.  In the
case where $X$ represents the germ of the whole affine space
$\IC^a\x\IC^b$ and $Y$ is the critical set $\Sigma(f)$, L\^e and Saito
\cite{LS73} proved that if the number of vanishing cycles, the Milnor
number $\mu(Z_y)$ at $0$, is constant in $y$, then \Af\ holds.

We prove a definitive generalization of this celebrated theorem of L\^e
and Saito to ICIS germs.  (However, Green and Massey, \cite{GnM96} and
\cite{M98}, showed that other information about the vanishing cycles
implies \Af\ for families with generalized isolated singularities.)
Namely, we characterize \Af\ by the constancy in $y$ of the Milnor
numbers $\mu(X_y)$ and $\mu(Z_y)$.  Notice that these Milnor numbers
refer to the ambient topology of $X_y$ and $Z_y$ in $\IC^a$, but the
isotopy theorem does not.

We characterize \Wf\ similarly, by the constancy of the sequences of
Milnor numbers $\mu_i(X_y)$ and $\mu_i(Z_y)$ of the sections of $X_y$
and $Z_y$ by general linear spaces of codimension $i$ for
$i=0,\dots,a-k$.  The appearance of these linear sections is a
reflection of the Lipschitz-like nature of the vector fields that we
integrate to prove the isotopy theorem.  (The corresponding theorem for
W was proved for hypersurface germs in part by Teissier and in part by
Brian\c con and Speder, and for general ICIS germs by Gaffney
\cite{Ga96, Thm.~1}.)

We also characterize \Af\ and \Wf\ by the constancy of certain
Buchsbaum--Rim multiplicities.  Say $X$ is defined in $\IC^a\x\IC^b$ by
the vanishing of $f_1,\dots,f_k$, and form the Jacobian matrix of
$f_1,\dots,f_k,f$ with respect to the first $a$ variables.  For each
$y$, the $a$ columns of this matrix generate a module over the local
ring $\cO_{X_y,0}$.  This is a submodule $\cM_y$ of finite colength of
the free module of rank $k+1$.  We characterize \Af\ by the constancy in
$y$ of the Buchsbaum--Rim multiplicity $e(\cM_y)$.  Let $\Im_y$ be the
maximal ideal of $X_y$; its appearance is a more direct reflection of
the Lipschitz-like nature of the vector fields.  We characterize \Wf\ by
the constancy of the Buchsbaum--Rim multiplicity of the product,
$e(\Im_y\cM_y)$.

Thus our theorems provide necessary and sufficient conditions for \Af\
(resp., for \Wf) to hold.  These conditions require the constancy in $y$
of certain numerical invariants of the individual pairs $X_y,f_y$.
Therefore, when these invariants are constant, then, in whatever way
the individual pairs are glued together to form a pair $X$, $f$,
necessarily \Af\ (resp., \Wf) holds.  Conversely, if one pair $X$, $f$
exists for which \Af\ (resp., \Wf) holds, then these invariants of
$X_y$, $f_y$ have the same values for all $y\in Y$ near 0.  Thus these
numerical invariants may be considered as indicators of
``\Af-equisingularity'' (resp., ``\Wf-equisingularity'').

The invariants depend only on the germs of $X_y$ and $f_y$ at $0$.  Yet,
of course, there is a significant difference between these germs and the
larger representatives themselves.  Assume that 0 is the only critical
point of the germ of $f_y$.  Still, the representative $f_y$ may have an
additional critical point no matter how close $y$ is to $0$.  However,
$f_y$ has no additional critical point if our numerical invariants are
constant.  Furthermore, if there is no such critical point, then \Af\
holds.  Conversely, if \Af\ holds, then, for all $y$ close to $0$, each
additional critical point of $f_y$ is also a critical point of $f$.

Our treatments of \Af\ and \Wf\ run in parallel.  Not only are the
statements similar, but, to a fair extent, the proofs are similar.
First, we prove that the two numerical characterizations of \Af\ (resp.,
of \Wf) are equivalent by proving that the constancy of the Milnor
numbers is equivalent to the constancy of the Buchsbaum--Rim
multiplicities.  This proof involves the theorem of L\^e and Greuel,
which re-expresses the Milnor numbers algebraically, and some theorems
of Buchsbaum and Rim, which re-express the multiplicity $e(\cM_y)$ as a
length (resp., and a more recent theorem, which re-expresses
$e(\Im_y\cM_y)$ as a linear combination of the polar multiplicities of
$X_y$).

We prove the characterization of \Af\ (resp., of \Wf) by Buchsbaum--Rim
multiplicities using the theory of integral dependence of elements on
submodules of free modules.  Let $\cM$ be the module, over the local
ring of $X$ at $0$, generated by the $a$ columns of the Jacobian matrix
above.  Let $g_j$ be the column vector of partial derivatives of
$f_1,\dots,f_k,f$ with respect to the $j$th coordinate variable on
$\IC^b$.  Finally, let $\Im_Y$ be the ideal of $Y$ in $X$.

We characterize \Af\ (resp., \Wf) by the integral dependence of
$g_1,\dots,g_b$ on $\cM$ (resp., on $\Im_Y\cM$).  We prove this
characterization of \Af\ using a remarkable result, Lemma (5.7) in
\cite{GM98}, concerning the geometry of the relative conormal variety.
Technically, this result is the new ingredient, which permits perfecting
the treatment of \Af\ in \cite{GK98}.  We prove the characterization of
\Wf\  via a direct computation with analytic inequalities.

The restrictions $g_j|X_y$ are always integrally dependent on $\cM_y$
(resp., on $\Im_y\cM_y$) for all $y$ in a dense Zariski open subset of
$Y$, as \Af\ (resp., \Wf) holds for all $y$ in such a subset by a
celebrated theorem proved by Hironaka (resp., by Henry, Merle, and
Sabbah).  It follows that, if the Buchsbaum-Rim multiplicity $e(\cM_y)$
(resp., $e(\Im_y\cM_y)$) is constant, then the $g_j$ are dependent on
$\cM$ (resp., on $\Im_Y\cM$) by virtue of the ``principle of
specialization of integral dependence.''

This principle is the main algebraic result, Theorem~(1.8), in
\cite{GK98}.  For ideals, it was discovered, named, and proved by
Teissier \cite{Te73, 3.2, p.~330} and \cite{Te80, App.~I}.  To prove it
for modules, we follow Teissier's approach, but first we must generalize
certain basic results in the theory of multiplicity for ideals,
including upper semicontinuity, Rees's characterization of integral
dependence, and B\"oger's generalization of it (Teissier rediscovered
the latter in the case at hand).  Upper semicontinuity of the
Buchsbaum--Rim multiplicity is proved in Proposition (1.1) of
\cite{GK98}.  The generalizations of Rees's theorem and of B\"oger's
theorem are proved in Theorems (6.7a)(iii) and (10.9) in \cite{KT94}.
Moreover, B\"oger's theorem is given a new proof in \cite{KT97}; see its
subsections (1.4) and (1.7).  This proof is simpler, shorter, and more
direct.  It was inspired by Lemma~(5.7) in \cite{GM98}.

The present article is meant to give a feeling for the nature and the
spirit of the work in \cite{GK98} and its extension in \cite{GM98}.  The
latter article also surveys a lot of other work in equisingularity
theory, and discusses its historical development.  This is a masterful
report, and may be highly recommended.

The present discussion is, of course, limited, and so achieves greater
focus.  Some results in \cite{GK98} are not discussed here, and some are
discussed in restricted generality.  Some proofs are abridged, and some
are omitted.  Also, the general theory is recalled in the special case
at hand.  However, some noteworthy results were not explicitly discussed
before, and they are stated and proved here.

In Section~2, we review the formal definitions of \Af\ and \Wf.  We
discuss isotopy, and state the two theorems, which characterize \Af\ and
\Wf\ by the constancy of Milnor numbers.  In Section~3, we review the
formal definition of the Buchsbaum--Rim multiplicity.  Then we
reformulate the theorems in terms of the constancy of Buchsbaum--Rim
multiplicities, and we prove that the two formulations are equivalent.
In Section~4, we review the theory of integral dependence of elements on
submodules of free modules, and we characterize \Af\ and \Wf\ in terms
of integral dependence.  Finally, we apply this theory to prove the
reformulated theorems.

\sct2.  Equisingularity

 Let $X$ be a complex-analytic germ at $0$ in $\IC^a\x\IC^b$.  Say that,
on a (Euclidean) neighborhood of $0$, we have
	$$X:f_1,\dots,f_k=0,$$
 where each $f_i$ is an analytic function $f_i(x,y)$ of the two sets of
variables,
        $$x=(x_1,\dots,x_a)\and y=(y_1,\dots,y_b).$$
 Assume that $X$ is a reduced complete intersection of codimension $k$
with $k<a$.

For fixed $y$, let $X_y\subset\IC^a$ denote the locus of $x$ such that
$(x,y)\in X$; so $X_y$ is the locus of zeros of the $f_i(x,y)$ as $x$
varies.  Assume that, if $X_y$ is nonempty, then $0\in X_y$.  Let $Y$ be
the locus of $y$ with $(0,y)\in X$, assume that $Y$ contains
 a neighborhood of $0$ in $\IC^b$, and identify $Y$ with $0\x Y$.  View
$Y$ as the parameter space and $X$ as the total space of the family of
$X_y$.  Assume that the $X_y$ represent germs at $0$ of isolated
complete-intersection singularities (ICIS germs) of codimension $k$.

Let $f$ be a nonconstant analytic function on $X$.  Set
 $$\textstyle f_y:=f|X_y\and Z:=X\cap f^{-1}0\and Z_y:=X_y\cap f^{-1}0.$$
 Assume that $f$ vanishes on $Y$, so $Y\subset Z$.

If $x\in X$ is a simple point of the level hypersurface $f^{-1}fx$, then
$x$ must be a simple point of $X$.  Let $\Sigma(f)$ denote the critical
set, the union of the singular sets of the various level hypersurfaces.
Then, in other words, $\Sigma(f)$ contains the singular set of $X$.

In turn, $\Sigma(f)$ is contained in the union of the critical sets of
the various restrictions $f_y$.  Denote this union by $\Sigma_Y(f)$.
Note that $\Sigma_Y(f)$ consists of all the singular points of all the
level hypersurfaces of all the restrictions $f_y$ for all the $y$ in
$Y$.  Replacing $X$ by a smaller representative of the same germ, we may
also assume that every component of $\Sigma_Y(f)$ contains 0.
Similarly, if we assume that $Z_0$ represents a germ with an isolated
singularity at $0$, then we may also assume that the projection
$\Sigma_Y(f)\to Y$ is finite.  Note that we work only with the reduced
structures on the sets $\Sigma(f)$ and $\Sigma_Y(f)$.

The {\it Thom condition\/} \Af\ can be formulated succinctly using the
{\it relative conormal variety\/} $C(X,f)$ and the {\it absolute
conormal variety} $C(Y)$.  These varieties are defined as follows.  Both
are closures in $X\x\IP^{a+b-1}$, or more correctly, in the restriction
to $X$ of the projectivized cotangent bundle of $\IC^a\x\IC^b$.  The
former closure is that of the set of pairs $(x,H)$ such that $x$ is a
point in $X-\Sigma(f)$ and $H$ is a hyperplane tangent at $x$ to the
level surface $f^{-1}fx$.  The latter closure is that of the set of
pairs $(x,H)$ such that $x\in Y\subset H$.  In these terms, \Af\ is said
to be satisfied by the pair $(X-\Sigma(f),Y)$ at 0 if the fiber of
$C(X,f)$ over $0\in X$ lies in $C(Y)$.

Suppose the germs of $\Sigma(f)$ and $Y$ at 0 are equal.  If \Af\ holds
at 0, then  \Af\ holds at every $y$ in a neighborhood $U$ of 0 in
$Y$.  This statement is not obvious from the definition, but follows
from Proposition~(4.2)(2) below.

Hence, $C(Y)$ contains the preimage of $U$ in $C(X,f)$, which will be
denoted by $C(X,f)|U$.  Since the intersection $C(X,f)\cap C(Y)$ always
projects into $Y$, we conclude that \Af\ holds at 0 if and only if,
after we replace $X$ by a smaller representative of the same germ, we
obtain the set-theoretic equation,
	$$C(X,f)\cap  C(Y)=C(X,f)|Y.$$
 In other words, {\it \Af\ holds at $0$ if and
only if, along the fiber of $C(X,f)$ over $0\in X$, the ideal of the
intersection $C(X,f)\cap C(Y)$ has the same radical as the ideal of the
preimage $C(X,f)|Y$}; compare with the {\it Remarque\/} on \p.550 in
\cite{LT88}.

The condition \Af\ can also be expressed analytically in terms of the
``angular distance'' $\dist(Y,T_xf^{-1}fx)$ from $Y$ to the
tangent space at $x$ to the level hypersurface.  Namely, \Af\ holds at 0
if and only if this distance approaches 0 as $x$ approaches 0 along any
analytic path $\phi\:(\IC,0)\to(X,0)$ such that $\phi(u)$ lies in
$X-\Sigma(f)$ for $u\neq0$.  Now, this distance approaches 0 if and only
if the inequality,
	$$\dist(Y, T_xf^{-1}fx)\le c\cdot\dist(x,Y)^e,$$
 holds where the constant $c$ and the exponent $e$ depend on the path
$\phi$.

The {\it Whitney condition relative to\/} $f$, denote \Wf, is defined by
requiring the preceding analytic inequality to hold for some constant
$c$ independent of $\phi$ and for $e=1$.  This condition generalizes
Teissier's condition of `c-equisingularity' (see \cite{LT88, top,
p.~550}).  It reduces to Whitney's Condition B, in the equivalent form
of Verdier's condition $W$ \cite{V76, Sect.~1}, when $f$ is constant
(and so vanishes).

There is also a characterization of \Wf\ in terms of the conormal
varieties $C(X,f)$ and $C(Y)$, which strengthens the corresponding
characterization of \Af\ in an interesting way.  Indeed, L\^e and
Teissier proved (a more general version of) the following result in
their Proposition~1.3.8 on \p.550 of \cite{LT88} (see Proposition~(6.1)
in \cite{GK98} for another treatment of the general case): {\it \Wf\
holds at $0$ if and only if, along the fiber of $C(X,f)$ over $0\in X$,
the ideal of the intersection $C(X,f)\cap C(Y)$ is integral over the
ideal of the preimage $C(X,f)|Y$.}

The condition \Wf\ implies that the pair $X,f$ is topologically right
trivial over $Y$.  Indeed, Thom, Mather, and Teissier, Verdier, Gaffney
and others introduced and developed the requisite methods to prove this
triviality by integrating vector fields.  Namely, it is possible to take
the constant tangent vector field to $Y$ and carefully lift it to $X$ so
that the lift is Lipschitz-like (Fr.~{\it rugueux}), hence integrable,
and is tangent to the level hypersurfaces of $f$, so that the integral
gives an appropriate continuous flow on $X$.

  Thus we obtain the following isotopy theorem: {\it If the critical set
$\Sigma(f)$ is equal to $Y$ and if the pair $(X-Y,Y)$ satisfies \Wf\ at
$0$, then, after $X$ is replaced by a neighborhood of $0$, there is a
homeomorphism $h\:X_0\x Y\to X$ such that $fh=f_0\x1_Y$ and $h$ is
$C^\infty$ off Y.}

Let $X^i$ be the section of $X$ by a general linear space of codimension
$i$ containing $Y$, and set $f^i:=f|X^i$; so $X^0=X$ and $f^0=f$.  If
\Wf\ holds, then \Wfi\ holds for $i=0,\dots,a-k-1$ because the requisite
analytic inequality follows from that for \Wf.  Hence the pair
$X^i,\,f^i$ is also topologically trivial by the isotopy theorem.

The next theorem characterizes \Wf\ in terms of the Milnor numbers
$\mu_i(X_y)$ and $\mu_i(Z_y)$ of the sections of $X_y$ and $Z_y$ by a
general linear space through 0 of codimension $i$ for $i=0,\dots,a-k$.
By convention, $\mu_{a-k}(X_y)$ and $\mu_{a-k-1}(Z_y)$ are the ordinary
multiplicities at $0$ diminished by $1$, and $\mu_{a-k}(Z_y)=1$.
However, for $\mu_i(Z_y)$ to be defined, $Z_y$ too must have an isolated
singularity at $0$.  For $y$ near $0$, it does if the germs of
$\Sigma(f)$ and $Y$ at $0$ are equal, and if \Wf, or simply \Af, holds;
indeed, then the germs of $\Sigma(f)$ and $\Sigma_Y(f)$ at $0$ are equal
by Lemma~(4.3) below.

 \thm1 The following three conditions are equivalent:\smallbreak
 \item i the germs of\/ $\Sigma(f)$ and $Y$ at $0$ are equal, and the
pair $(X-\Sigma(f),Y)$ satisfies \Wf\ at $0;$
 \item ii for the $y$ in a  neighborhood of\/ $0$ in $Y$, the
level hypersurface $Z_y$ has an isolated singularity at $0$, and the
sequences of Milnor numbers, $\{\mu_i(X_y)\}$ and $\{\mu_i(Z_y)\}$, are
constant in $y;$
 \item ii$'$ for the $y$ in a  neighborhood of\/ $0$ in $Y$, the
level hypersurface $Z_y$ has an isolated singularity at $0$, and the
sequence of sums of Milnor numbers, $\{\mu_i(X_y)+\mu_i(Z_y)\}$, is
constant in $y$.
 \endproclaim

This theorem is part of Theorem~(6.4) in \cite{GK98}.  The latter also
asserts that (i) and (ii) are equivalent to the following condition:
{\it the germs of $\Sigma_Y(f)$ and $Y$ at $0$ are equal, and both pairs
$(X-Y,Y)$ and $(Z-Y,Y)$ satisfy the absolute Whitney condition W at
$0$.}

This additional condition is implied by (i).  Indeed, the two germs are
equal by Lemma~(4.3) below because \Wf\ implies \Af.   Furthermore,
 (i) implies the
requisite analytic inequalities.  Indeed, $T_xf^{-1}fx\subset
T_xX$ and, if $x\in Z$, then $T_xf^{-1}fx=T_xZ$.

The additional condition implies (ii); indeed, this implication is
virtually the assertion of Th\'eor\`eme~(10.1) on \p.223 of \cite{N80}.
Trivially,  (ii) implies (ii$'$).

Finally, (ii$'$) implies (i), but our proof is more involved, and runs
as follows: in Theorem~(3.1) below, we replace (ii$'$) with an
equivalent condition involving a Buchsbaum--multiplicity, and then
at the end of Section~4, we prove that latter implies (i).  On the other
hand, there is a direct proof that the additional condition implies (i);
indeed, Brian\c con, Maisonobe and Merle gave such a proof in \cite{BMM,
Thm.~4.3.2, p.~543} in a more general setting using a different
approach.

The next theorem characterizes \Af\ in terms of Milnor numbers.

\thm2 The following four conditions are equivalent:
 \smallbreak
 \item i the germs of\/ $\Sigma(f)$ and $Y$ at $0$ are equal, and the pair
$(X-\Sigma(f),Y)$ satisfies \Af\ at $0;$
 \item ii the germs of\/ $\Sigma_Y(f)$ and $Y$ at $0$ are equal;
 \item iii for the $y$ in a  neighborhood of $0$ in $Y$, the
level hypersurface $Z_y$ has an isolated singularity at $0$, and the
Milnor numbers, $\mu(X_y)$ and $\mu(Z_y)$, are constant in $y;$
 \item iii$'$ for the $y$ in a  neighborhood of $0$ in $Y$, the
level hypersurface $Z_y$ has an isolated singularity at $0$, and the sum
of Milnor numbers, $\mu(X_y)+\mu(Z_y)$, is constant in $y$.
 \endproclaim

Notice that (iii) implies (iii$'$) trivially.  The converse holds
because $\mu(X_y)$ and $\mu(Z_y)$ are each upper semicontinuous by
\cite{Lo84, bot.~p.~126}.  We prove the rest of the theorem in two steps
too: first, we replace (iii$'$) with an equivalent condition involving a
Buchsbaum--Rim multiplicity, obtaining Theorem~(3.2); then at the end of
Section~4, we prove Theorem~(3.2).

In the case where where $X$ represents the germ of $\IC^a\x\IC^b$, L\^e
and Saito \cite{LS73} proved that (iii) implies (i).  They used Morse
theory, but Teissier reproved their theorem almost immediately using
more algebraic-geometric methods.  Teissier's work served as a model for
most of the work  in this report.

Combining the preceding two theorems, we obtain the following corollary,
which asserts the equivalence of \Wf,  \Wfi, and \Afi.

 \cor3 If $\Sigma(f)=Y$, then the following conditions are equivalent:
 \smallbreak
 \item i the pair $(X-Y,Y)$ satisfies \Wf\ at $0;$
 \item ii the pair $(X^i-Y,Y)$ satisfies \Wfi\ at $0$ for
every $i;$
 \item iii the pair $(X^i-Y,Y)$ satisfies  \Afi\ at $0$ for
every $i$.
 \endproclaim

\sct3. Multiplicity

In this section, we reformulate Theorems (2.1) and (2.2).  Instead of
Milnor numbers, we use certain Buchsbaum--Rim multiplicities.  Thus we
obtain Theorems (3.1) and (3.2).  In the next section, we discuss the
proofs of these reformulated theorems.

The Buchsbaum--Rim multiplicity was introduced by Buchsbaum and Rim
\cite{BR63} in 1963, and the theory has been developed more recently by
Kirby and Rees \cite{KR94}, by Henry and Merle \cite{HM94}, and by
Thorup and the author \cite{KT94} and \cite{KT96}.  For our purposes
here, it suffices to work over the local ring $\cO$ of a
complex-analytic germ, say one of dimension $d$.  Let $\cE$ be a free
$\cO$-module, and $\cM$ a submodule of finite colength; that is, the
vector-space dimension $\dim_{\IC}(\cE/\cM)$ is finite.

The Buchsbaum--Rim multiplicity generalizes the ordinary multiplicity.
In the case where $\cE$ is the ring $\cO$ and where $\cM$ is an ideal
$\cI$, Samuel defined the multiplicity $e(\cI)$ in 1951 to be the
rectified leading coefficient $e$ of the Hilbert--Samuel polynomial,
	$$\dim_{\IC}(\cE/\cI^n)=e\,n^d/d!+\cdots
		\hbox{ for }n\gg0.$$

Buchsbaum and Rim considered the case in which $\cE$ is free of
arbitrary (finite) rank $r$.  They generalized Samuel's definition
essentially as follows.  Form the symmetric algebra $\cO[\cE]$; it is
just the polynomial algebra in $r$ variables with coefficients in $\cO$.
Form the {\it Rees\/} algebra $\cO[\cM]$; it is just the subalgebra
generated by $\cM$ placed in degree 1.  Both these algebras are graded;
denote their $n$th graded pieces by $\cO[\cE]_n$ and $\cO[\cM]_n$.  For
example, if $\cE=\cO$ and $\cM=\cI$, then $\cO[\cE]_n=\cO$ and
$\cO[\cM]_n=\cI^n$.

Buchsbaum and Rim formed the quotient of these two graded pieces, and
proved that its dimension is eventually given by a polynomial in $n$ of
degree $d+r-1$; thus,
	$$\dim_{\IC}(\cO[\cE]_n/\cO[\cM]_n)=e\,n^{d+r-1}/(d+r-1)!+\cdots
		\hbox{ for }n\gg0.$$
 Then they defined the {\it multiplicity\/} $e(\cM)$ to be $e$.

To use the Buchsbaum--Rim multiplicity, return to the setup described at
the beginning of Section~2.  Fix an extension of $f$ over a neighborhood
of $X$ in $\IC^a\x\IC^b$ on which $f_1,\dots,f_k$ are defined, and
abusing notation, denote the extension too by $f$.  Form the Jacobian
matrix with respect to the first set of variables:
	$$\left[\matrix{\pd f_1/x_1 &\ldots&\pd f_1/x_a \cr
			\vdots	    &\ddots&\vdots\cr
			\pd f_k/x_1 &\ldots&\pd f_k/x_a \cr
			\pd f/x_1   &\ldots&\pd f/x_a \cr}\right].$$
 Its columns generate a module over the affine ring $\cO_{X}$,
	$$JM_x(f_1,\dots,f_k;f)\subset\cE:=\cO_{X}^{k+1},$$
 called the {\it Jacobian module\/} of $f_1,\dots,f_k;f$ with
respect to the $x$-variables.

Given $y\in Y$, form the image $\cM_y$ of $JM_x(f_1,\dots,f_k;f)$ in the
free module $\cE_y$ of rank $k+1$ over the local ring $\cO_{X_y,0}$; so
	$$\cM_y:=JM_x(f_1,\dots,f_k;f)|X_y\subset\cE_y.$$
 In addition, denote the maximal ideal of $\cO_{X_y,0}$ by $\Im_y$.

Suppose for a moment that each $Z_y$ has an isolated singularity at
$0$. Since
	$$\Sigma_Y(f)=\Supp\bigl(\cE/JM_x(f_1,\dots,f_k;f)\bigr),$$
the module $\cM_y$ has finite colength.  Hence the Buchsbaum--Rim
multiplicities $e(\cM_y)$ and $e(\Im_y\cM_y)$ are defined.

We can now reformulate Theorems (2.1) and (2.2); we obtain
Theorems~(3.1) and (3.2).  After stating them, we discuss the proof that
the two formulations are equivalent.

 \thm1 The following two conditions are equivalent:\smallbreak
 \item i the germs of\/ $\Sigma(f)$ and $Y$ at $0$ are equal, and the pair
$(X-\Sigma(f),Y)$ satisfies \Wf\ at $0;$
 \item ii for the $y$ in a  neighborhood of $0$ in $Y$, the
level hypersurface $Z_y$ has an isolated singularity at $0$, and the
multiplicity $e(\Im_y\cM_y)$ is constant in $y$.
 \endproclaim
\thm2 The following three conditions are equivalent: \smallbreak
 \item i the germs of\/ $\Sigma(f)$ and $Y$ at $0$ are equal, and the pair
$(X-\Sigma(f),Y)$ satisfies \Af\ at $0;$
 \item ii the germs of\/ $\Sigma_Y(f)$ and $Y$ at $0$ are equal;
 \item iii for the $y$ in a  neighborhood of $0$ in $Y$, the
level hypersurface $Z_y$ has an isolated singularity at $0$, and the
multiplicity $e(\cM_y)$ is constant in $y$.
 \endproclaim

The new theorems are equivalent to the old because of the following
lemma.  Indeed, the summands in (3.3.2) are upper semicontinuous; so
they are constant if and only if $e(\Im_y\cM_y)$ is.  For future use,
note that therefore the lemma yields this: {\it if $e(\Im_y\cM_y)$ is
constant, then so is $e(\cM_y)$.}

\lem3 For each $y\in Y$, we have these two equations:
    $$\displaylines{e(\cM_y)=\mu(X_y)+\mu(Z_y);\dleqno{(3.3.1)}\cr
e(\Im_y\cM_y)=\sum_{i=0}^{a-k}{\textstyle{a-1\choose i}}
       \bigl(\mu_i(X_y)+\mu_i(Z_y)\bigr). \dleqno{(3.3.2)}\cr}$$
 Moreover, each sum $\mu_i(X_y)+\mu_i(Z_y)$ is upper semicontinuous in
$y$.
 \endproclaim

To prove Equation (3.3.1), note that $X_y$ is a complete intersection of
dimension $a-k$ and that $\cM_y$ is a submodule of $\cO_{X_y,0}^{k+1}$
generated by $a$ elements, namely, the columns of the Jacobian matrix.
Denote the ideal of maximal minors of this matrix by $\cJ$.  Then some
theorems of Buchsbaum and Rim \cite{BR63, 2.4, 4.3, 4.5} yield
$$e(\cM_y)=\dim_\IC\bigl(\cO_{X_y,0}/\cJ\cO_{X_y,0}\bigr).\eqno(3.3.3)$$
 The right side is equal to the sum $\mu(X_y)+\mu(Z_y)$ by the theorem
of L\^e \cite{Le74, Thm.~3.7.1, p.~130} and Greuel \cite{Gr75, Kor.~5.5,
p.~263}.  Thus Equation (3.3.1) holds.

To prove Equation (3.3.2), for each $i$, let $P_i$ be a general linear
space of codimension $i$ in $\IC^a$, and let $\Pi^i$ be the
$i$-dimensional ``relative polar subscheme'' of $f_y$ with $P_i$ as
pole.  By definition, $\Pi^i$ is the closure in $X_y$ of the locus of
simple points $x$ of the level hypersurface surface $X_y\cap f^{-1}fx$
such that there exists a tangent hyperplane at $x$ that contains $P_i$.
Algebraically, $\Pi^i$ is cut out of $X_y$ by the maximal minors of the
Jacobian matrix of the map $\IC^a\to\IC^{k+1}\x\IC^i$ with components
$f_1,\dots,f_k,f$ and $p_i$, where $p_i\:\IC^a\to\IC^i$ is the linear
map with kernel $P_i$.  Hence $\Pi^i$ is Cohen--Macaulay as $X_y$ is.

{}From the polar multiplicity formula \cite{KT94, Thm.~(9.8)(i)} (compare
with \cite{HM94, 4.2.7} and \cite{Ga96, \S3}), it follows that
	$$e(\Im_y\cM_y)= \sum_{i=0}^{a-k}\tc a-1 i m(\Pi^i),$$
 where $m(\Pi^i)$ is the ordinary multiplicity at $0$ of $\Pi^i$.  In
particular, $m(\Pi^{a-k})$ is simply the multiplicity of $X_y$ at $0;$
so it is equal to $\mu_{a-k}(X_y)+1$.  For any $i$, since $\Pi^i$ is
Cohen--Macaulay,
	$$m(\Pi^i) = \dim_\IC(\cO_{\Pi^i,0}/\cI_i),\dno3.4$$
 where $\cI_i$ is the ideal of any linear space $L_i$ of codimension $i$
in $\IC^a$ that is transverse to $\Pi^i$.  It follows that $m(\Pi^i)$ is
upper semicontinuous.

Remarkably, although $\Pi^i$ is defined using $P_i$, nevertheless $P_i$
is transverse to $\Pi^i$ simply because $P_i$ is general.  This
important result was proved by Teissier in \cite{Te81, (4.1.8), p.~569}
as a consequence of his general idealistic Bertini theorem.  The result
was also proved, at about the same time, by Henry and Merle \cite{HM82,
Cor.~2, p.~195}.  The result is reproved in Lemma~(6.2) of \cite{GK98}
in a new way, using the theory of the \Wf\ condition in the spirit of
the current work.  By this transversality result, we may take $L_i$ to
be $P_i$.  Then, for $i<a-k$, the right side of Equation~(3.3.4) is
equal to $\mu_i(X_y)+ \mu_i(Z_y)$ by the theorem of L\^e and Greuel.
The asserted formula follows immediately, and the proof is complete.

\sct4. Dependence

In this section, we discuss the proofs of Theorems~(3.1) and (3.2), and
thereby complete the proofs of Theorems (2.1) and (2.2).  Our main
technical tool is the theory of integral dependence of elements on
modules, which generalizes the older theory for ideals.  We begin by
reviewing this theory.

Let $\cO$ be the local ring of a complex-analytic germ $(X,0)$.  Let
$\cE$ be a free $\cO$-module, $\cM$ a submodule, and $g\in\cE$ an
element.  By definition, $g$ is {\it integrally dependent\/} on $\cM$
if, when $g$ is viewed as an element of degree $1$ in the symmetric
algebra $\cO[\cE]$, then $g$ is integrally dependent on the Rees algebra
$\cO[\cM]$ in the usual sense; namely, $g$ satisfies an equation of
integral dependence,
	$$g^n+r_1g^{n-1}+\cdots+r_n=0,$$
 where $n\ge1$ and $r_i\in\cO[\cM]$, both depending on $g$.  Of course,
each $r_i$ may be replaced by its homogeneous piece of degree $i$ if
desired.

There are two useful criteria for integral dependence.  The first is a
form of the valuative criterion, but it is also known in the trade as
the ``curve criterion.'' \smallskip
 (Curve criterion)\enspace {\it For $g\in\cE$ to be integrally
dependent on $\cM$  it is necessary that, for every map germ
$\phi\:(\IC,0)\to(X,0)$, the pullback $\phi^*g$ lie in the pullback
$\phi^*\cM$, viewed in the free $\cO_{\IC,0}$-module $\phi^*\cE$, or
put more concisely,
		 $$\phi^*g\in\phi^*\cM\subset\phi^*\cE;$$
 conversely, it is sufficient that this condition obtain for every
nonconstant $\phi$ whose image meets any given dense Zariski open subset
of $X$.}
 \smallskip
The second criterion is an analytic inequality, which re-expresses
integral dependence in terms of speeds of vanishing. \smallskip
 (Analytic criterion)\enspace {\it For $g\in\cE$ to be integrally
dependent on $\cM$ it is necessary that, for any finite set of generators
$g_i$ of $\cM$, there exist a Euclidean neighborhood $U$ of\/ $0$ in $X$
and a constant $c$ such that $|g(x)|\le c\,\max|g_i(x)|$ for any $x\in
U$; conversely, it is sufficient that this condition obtain for some
finite set of generators $g_i$ of $\cM$.}
 \smallbreak Indeed, the condition in the ``curve criterion'' is
equivalent to the condition in the ``analytic criterion'' by
Proposition~1.11 on \p.306 of \cite{Ga92}.  In fact, on \p.303 in
\cite{Ga92}, Gaffney takes the former as the defining condition of
integral dependence.  On \p.305, he proves that this definition is
equivalent to Rees's definition \cite{Re87, \p.435}.  Finally,
Theorem~1$\cdot$5 in \cite{Re87, \p.437} yields the curve criterion.

In fact, neither \cite{Ga92} nor \cite{Re87} treat the present version
of the curve criterion, with its weaker sufficiency condition.  However,
it is not difficult to extend that work: if $\phi^*g\notin\phi^*\cM$,
then $\phi$ can be tweaked, preserving this relation, so that its image
does meet the given open set (see the proof of Prop.~1.7 on \p.304 in
\cite{Ga92}).

 \smallbreak {\bf From now on,} work in the setup described at the
beginning of Section~2.  The next lemma is the main algebraic result,
Theorem~(1.8), in \cite{GK98}.  Let $\cE$, $\cM$, and $g$ be as above,
but assume that they arise from a free module, a submodule, and an
element defined over the affine ring of $X$.  In addition, for each
$y\in Y$, form the restriction $\cE_y$ of the free module, the image
$\cM_y\subset \cE_y$ of the submodule, and the image $g_y\in\cE_y$ of
the element.  Finally, assume that $\cM_y$ has finite colength.

\lem1 \(Principle of specialization of integral dependence)\enspace
 Assume that the Buchsbaum--Rim multiplicity $e(\cM_y)$ is constant in
$y$.  Then $g$ is integrally dependent on $\cM$ if $g_y$ is integrally
dependent on $\cM_y$ for all $y$ in a dense Zariski open subset of $Y$.
 \endproclaim

Of course, if $g$ is integrally dependent on $\cM$, then $g_y$ is
integrally dependent on $\cM_y$ for all $y$ in a  neighborhood of
$0$, but the latter condition is strictly weaker than the former.
Indeed, a simple  example is given in Example~(1.3) of
\cite{GK98}.  Thus the conclusion of the lemma is stronger than it might
seem at first.
 \smallskip
To prove Theorems~(3.1) and (3.2), take $\cM$ to be the Jacobian module,
	$$\cM:=JM_x(f_1,\dots,f_k;f)\subset\cE:=\cO_{X}^{k+1},$$
 the column space of the Jacobian matrix; see Section~3.  For
convenience, let $\cM$ also denote the induced module over the local
ring $\cO_{X,0}$.  For $1\le j\le b$, let $g_j$ be the column vector,
	$$g_j:=\left[\matrix{\pd f_1/y_j \cr
			\vdots	    \cr
			\pd f_k/y_j \cr
			\pd f/y_j \cr}\right].$$
 Finally, denote the ideal of $Y$ in $X$ by $\Im_Y$; so $\Im_Y:=
(x_1,\dots,x_a)\cO_X$.

The next result characterizes \Af\ (resp., \Wf) by the integral dependence
of the $g_j$ on the submodule  $\cM$ (resp.,  $\Im_Y\cM$) of $\cE$.

\prp2 Assume that $\Sigma(f)=Y$.
 \part1 Then $(X-Y,Y)$ satisfies \Wf\ at\/ $0$ if and only if the columns
$g_1,\dots,g_b$ are all integrally dependent on\/ $\Im_Y\cM$.
 \part2 Then $(X-Y,Y)$ satisfies \Af\ at\/ $0$ if and only if the columns
$g_1,\dots,g_b$ are all integrally dependent on $\cM$.
 \endproclaim

Indeed, (1) is part of Proposition~(6.1) in \cite{GK98}.  The
characterization is proved by developing the analytic inequalities
involved in the definition of \Wf\ until the condition in the analytic
criterion is met.

Part~(2) is not explicitly stated in either \cite{GK98} or \cite{GM98},
but may be derived as follows.  First note the formula (Gaffney,
priv.~comm., 1990),
	$$C(X,f):=\Projan(\cO_X[\cM,g_1,\cdots, g_b]).$$
 Indeed, both sides are closed subvarieties of $X\x\IP^{a+b-1}$.  Both
are equal, over $X-\Sigma(f)$, to the set of pairs $(x,H)$ such that $H$
is a hyperplane tangent at $x$ to the level surface $f^{-1}fx$.  The
left side is, by definition, the closure of this set.  The right side is
also the closure of this open subset of itself, because its algebra is,
by construction, a subalgebra of the symmetric algebra $\cO[\cE]$.

Let $V$ be an arbitrary component of the preimage of $Y$ in $C(X,f)$.
Since $\Sigma(f)=Y$, the dimension of $V$ is $a+b-1$ by Gaffney and
Massey's Lemma (5.7) in \cite{GM98}; see also Theorem~4.2 in \cite{M98}
and the corollary in \cite{KT97, (1.2)}.

Assume that $g_1,\cdots, g_b$ are integrally dependent on $\cM$.  Then,
after $X$ is replaced by a neighborhood of $0$, the inclusion of
$\cO_X[\cM]$ into $\cO_X[\cM,g_1,\cdots, g_b]$ induces a finite map,
	$$\gamma\:C(X,f)\to X\x\IP^{a-1},$$
 since $\cM$ is generated by $a$ elements.  Hence $\dim\gamma(V)=a+b-1$.
Since $\gamma(V)$ is contained in $Y\x\IP^{a-1}$, therefore these two
sets are equal.  Hence $V$ maps onto $Y$.  Since \Af\ holds at every $y$
in a dense Zariski open set $U$ of $Y$ by Hironaka's Theorem~2 on \p.247
in \cite{H76}, the preimage $V|U$ lies in $C(Y)$.  Since $V|U$ is
nonempty and so dense in $V$, therefore $V\subset C(Y)$.  Thus \Af\
holds at 0.

The converse is a special case of (1) of the following lemma, and
so the proof of the proposition is complete.

Both parts of the lemma are used below to complete the proofs of
Theorems~(3.1) and (3.2).  Moreover, both are interesting in their own
right.

\lem3 \(1) If  $(X-\Sigma(f),Y)$ satisfies \Af\ at\/ $0$, then
each $g_j$ is integrally dependent on $\cM$.
 \part2 If each $g_j$ is integrally dependent on $\cM$, then the germs
of\/ $\Sigma(f)$ and\/ $\Sigma_Y(f)$ at\/ $0$ are equal.
 \endproclaim

To prove (2), suppose that the germ of $\Sigma_Y(f)$ is strictly larger
than that of $\Sigma(f)$.  Then there is a path $\phi\:(\IC,0)\to(X,0)$
whose image lies in the former set, but outside the latter.  Now,
	$$\Sigma_Y(f)=\Supp(\cE/\cM)
	     \and\Sigma(f)=\Supp\bigl(\cE/(\cM,g_1,\dots,g_b)\bigr).$$
 Hence, for some $j$, the pullback $\phi^*g_j$ does not lie in the
pullback $\phi^*\cM$.  So this $g_j$ is not integrally dependent on
$\cM$ by the curve criterion.  Thus (2) holds.

To prove (1), take a $\phi\:(\IC,0)\to(X,0)$ whose image lies outside
$\Sigma_Y(f)$, so outside $\Sigma(f)$.  The gradients of
$f_1,\dots,f_k,f$ define hyperplanes tangent to the level hypersurfaces
of $f$.  Each hyperplane must approach along $\phi$ a hyperplane that
contains $Y$ since \Af\ holds.  Consider the last $b$ components of each
gradient.  Each such component must, therefore, vanish at $0$ along $\phi$
to order higher than the order of one, or more, of the first $a$
components.  Denote the maximal ideal of $(\IC,0)$ by $\Im$.  Then,
therefore,
	$$\phi^*g_1,\dots,\phi^*g_b\in\Im\cM\subset\cM.$$
  So, by the curve criterion, $g_j$ is integrally dependent on $\cM$.
Thus (1) holds.
 \medbreak
 To prove Theorem~(3.2), let $\cJ$ be the ideal of maximal minors of the
Jacobian matrix in Section~3.  Then
	$$\Sigma_Y(f)=\Supp(\cE/\cM)=\Supp(\cO_X/\cJ).$$
 Now, if either (ii) or (iii) holds, then each $Z_y$ has an isolated
singularity; hence, we may assume that $\Sigma_Y(f)$ is finite over $Y$.
Consequently, $\cO_X/\cJ$ is determinantal, so a Cohen--Macaulay ring by
Eagon's theorem.  Therefore, $\cO_X/\cJ$ is a Cohen--Macaulay
$\cO_Y$-module, and so, by the Auslander--Buchsbaum formula, free.
Hence the following number is constant in $y$:
	$$e'(y):=\dim_\IC(\cO_X/\cJ)(y).$$

For each $y$, this number $e'(y)$ is a sum of positive numbers, one for
each point $z$ in $\Sigma_Y(f)$ lying over $y$, and the number
corresponding to $z=0$ is equal to $e(\cM_y)$ by Equation~(3.3.3).
Hence, we have
	$$e(\cM_0)=e'(0)=e'(y)\ge e(\cM_y),$$
 with equality at the end for every $y$ if and only if $\Sigma_Y(f)=Y$.
Therefore, (ii) and (iii) are equivalent.

Assume (i).  Then (ii) follows directly from Lemma~(4.3).

Conversely, assume (ii).  Then, since $\Sigma_Y(f)$ and $\Sigma(f)$ and
$Y$ are nested, all three represent the same germ.  Furthermore, (iii)
holds; so $e(\cM_y)$ is constant.  Now, to prove that \Af\ holds, we use
Proposition~(4.2)(2).  We have to show that the $g_j$ are integrally
dependent on $\cM$.  However, they are so by the principle of
specialization of integral dependence.  Indeed, $(X-\Sigma(f),Y)$
satisfies \Af\ at every $y$ in a dense Zariski open subset $U$ of $Y$ by
Hironaka's Theorem~2 on \p.247 in \cite{H76}.  Hence, for $y\in U$, the
restrictions $g_j|X_y$ are integrally dependent on $\cM_y$ by
Proposition~(4.2)(2) again.  Thus Theorem~(3.2) is proved.

Finally, consider Theorem~(3.1).  We have left to prove that (ii)
implies (i).  So assume (ii); that is, $e(\Im_y\cM_y)$ is constant.
Then $e(\cM_y)$ is constant too; we noted this consequence before
stating Lemma~(3.3).  Hence, by Theorem~(3.2), the germs of $\Sigma(f)$
and $Y$ at $0$ are equal.  To prove that \Wf\ holds, we use
Proposition~(4.2)(1) and the principle of specialization of integral
dependence much as we just did for \Af; we need only replace Hironaka's
theorem by Henry, Merle, and Sabbah's Th\'eor\`eme~5.1 on \p.255 of
\cite{HMS84}.  (They attribute this result to Navarro, who didn't
publish it.)  Thus Theorem~(3.1) is proved.

\references

BMM
J. Brian\c con{,} P. Maisonobe and M. Merle,
 Localisation de syst\`emes diff\'erentiels, stratifications de Whitney
et condition de Thom,
 \invent 117 1994 531--550

BR63
 D.A. Buchsbaum and D.S. Rim,
 A generalized Koszul complex. II. Depth and multiplicity,
 \tams 111 1963 197--224

Ga92
 T. Gaffney,
 Integral closure of modules and Whitney equisingularity,
 \invent 107 1992 301--322

Ga96
 T. Gaffney,
 Multiplicities and equisingularity of ICIS germs,
\invent 123 1996 209--220

GK98
 T. Gaffney and S.L. Kleiman,
 Specialization of integral dependence for modules,
 {\it alg-geom/}9610003.

GM98
 T. Gaffney and D. Massey,
 Trends in equisingularity theory,
 in the proceedings of the Liverpool conference in honor of CTC Wall,
singularities volume, W. Bruce (ed.), Cambridge University Press, to
appear.

GnM96
 M.D. Green and D.B. Massey,
 Vanishing cycles and Thom's $a_f$ conditions, {\it Preprint\/} 1996.

Gr75
 G.M. Greuel,
 Der Gauss--Manin Zusammenhang isolierter Singularit\"aten von
voll\-st\"and\-igen Durchschnitten, {\it Dissertation}, G\"ottingen
(1973), \ma 214 1975 235--266

HM82
 J.P.G. Henry and M. Merle,
 Limites d'espaces tangents et transversalit\'e de vari\'et\'es
polaires,
 {\it in}  ``Proc. La R\'abida, 1981.'' J. M. Aroca, R. Buchweitz, M.
Giusti and M.  Merle (eds.) \splm 961 1982 189--199

HM94
 J.P.G. Henry and M. Merle,
 Conormal Space and Jacobian module.  A short dictionary,
 {\it in} ``Proceedings of the Lille Congress of Singularities,"
J.-P. Brasselet (ed.), London Math. Soc. Lecture Notes {\bf 201} (1994),
147--174.

HMS84
J.P.G. Henry, M. Merle,  and C. Sabbah,
 {\it Sur la condition de Thom stricte pour un morphisme analytique
complexe,}
 \asens 17 1984 227--268

H76
 H. Hironaka,
 Stratification and flatness,
 {\it in} ``Real and complex singularities, Nordic Summer School, Oslo,
1976,'' Sijthoff and Noordhoff, 1977, 199--265.

KR94
 D. Kirby and D. Rees,
{\it Multiplicities in graded rings I: The general theory,}
 in ``Commutative algebra: syzygies, multiplicities, and birational
algebra'' W. J. Heinzer, C. L. Huneke, J. D. Sally (eds.), \conm 159 1994
209--267

KT94
 S. Kleiman and A. Thorup,
 A geometric theory of the Buchsbaum--Rim multiplicity,
 \ja 167 1994 168--231

KT96
 S. Kleiman and A. Thorup,
 Mixed Buschsbaum--Rim Multiplicities,
 \ajm 118 1996 529--569

KT97
 S. Kleiman and A. Thorup,
 Conormal geometry of maximal minors,
 alg-geom/970818.

Le74
 D.T. L\^e,
 Calculation of Milnor number of isolated singularity of complete
intersection,
 \faa 8 1974 127--131

LS73
 D.T. L\^e and K. Saito,
 La constance du nombre de Milnor donne des bonnes stratifications,
 \crasp 277 1973 793--795

LT88
 D.T. L\^e and B. Teissier,
 Limites de'espaces tangent en g\'eom\'etrie analytique,
 \cmh 63 1988 540--578

Lo84
 E.J.N.  Looijenga,
 Isolated singular points on complete intersections,
 London Mathematical Society lecture note series {\bf 77},
Cambridge University Press, 1984.

M98
 D.B. Massey,
 Critical points of functions on singular spaces, {\it Preprint\/} 1998.

N80
 V. Navarro,
 Conditions de Whitney et sections planes,
 \invent 61 1980 199--226

Re87
 D. Rees,
 Reduction of modules,
 \mpcps 101 1987 431--449

Te73
 B. Teissier,
 Cycles \'evanescents, sections planes et conditions de Whitney,
 in ``Singularit\'es \`a Carg\`ese," Ast\'erisque {\bf 7--8} (1973),
285--362.

Te80
 B. Teissier,
 {\it R\'esolution simultan\'ee et cycles \'evanescents,}
 {\i in} ``S\'em. sur les singularit\'es des surfaces.''  Proc.
1976--77.  M.  Demazure, H. Pinkham and B. Teissier (eds.) \splm 777
1980 82--146

Te81
 B. Teissier,
 Multiplicit\'es polaires, sections planes, et conditions de Whitney,
 in ``Proc. La R\'abida, 1981.'' J. M. Aroca, R. Buchweitz, M. Giusti and
M.  Merle (eds.) \splm 961 1982 314--491

V76
 J.-L. Verdier,
 Stratifications de Whitney et th\'eor\`eme de Bertini--Sard,
\invent 36  1976 295--312

\endreferences
\mark{E}\bye